\magnification=1200

\catcode`\À=\active \defÀ{\`A}    \catcode`\à=\active \defà{\`a} 
\catcode`\Â=\active \defÂ{\^A}    \catcode`\â=\active \defâ{\^a} 
\catcode`\Æ=\active \defÆ{\AE}    \catcode`\æ=\active \defæ{\ae}
\catcode`\Ç=\active \defÇ{\c C}   \catcode`\ç=\active \defç{\c c}
\catcode`\È=\active \defÈ{\`E}    \catcode`\è=\active \defè{\`e} 
\catcode`\É=\active \defÉ{\'E}    \catcode`\é=\active \defé{\'e} 
\catcode`\Ê=\active \defÊ{\^E}    \catcode`\ê=\active \defê{\^e} 
\catcode`\Ë=\active \defË{\"E}    \catcode`\ë=\active \defë{\"e} 
\catcode`\Î=\active \defÎ{\^I}    \catcode`\î=\active \defî{\^\i}
\catcode`\Ï=\active \defÏ{\"I}    \catcode`\ï=\active \defï{\"\i}
\catcode`\Ô=\active \defÔ{\^O}    \catcode`\ô=\active \defô{\^o} 
\catcode`\Ù=\active \defÙ{\`U}    \catcode`\ù=\active \defù{\`u} 
\catcode`\Û=\active \defÛ{\^U}    \catcode`\û=\active \defû{\^u} 
\catcode`\Ü=\active \defÜ{\"U}    \catcode`\ü=\active \defü{\"u} 

\catcode`\ =\active \def { }
%         ^u            ^u^o 
% u = unbreakable space,  o = ordinary space

\hsize=11.25cm    
\vsize=18cm       
\parindent=12pt   \parskip=5pt     

\hoffset=.5cm   
\voffset=.8cm   

\pretolerance=500 \tolerance=1000  \brokenpenalty=5000

\catcode`\@=11

\font\eightrm=cmr8         \font\eighti=cmmi8
\font\eightsy=cmsy8        \font\eightbf=cmbx8
\font\eighttt=cmtt8        \font\eightit=cmti8
\font\eightsl=cmsl8        \font\sixrm=cmr6
\font\sixi=cmmi6           \font\sixsy=cmsy6
\font\sixbf=cmbx6

\font\tengoth=eufm10 
\font\eightgoth=eufm8  
\font\sevengoth=eufm7      
\font\sixgoth=eufm6        \font\fivegoth=eufm5

\skewchar\eighti='177 \skewchar\sixi='177
\skewchar\eightsy='60 \skewchar\sixsy='60

\newfam\gothfam           \newfam\bboardfam

\def\tenpoint{
  \textfont0=\tenrm \scriptfont0=\sevenrm \scriptscriptfont0=\fiverm
  \def\rm{\fam\z@\tenrm}
  \textfont1=\teni  \scriptfont1=\seveni  \scriptscriptfont1=\fivei
  \def\oldstyle{\fam\@ne\teni}\let\old=\oldstyle
  \textfont2=\tensy \scriptfont2=\sevensy \scriptscriptfont2=\fivesy
  \textfont\gothfam=\tengoth \scriptfont\gothfam=\sevengoth
  \scriptscriptfont\gothfam=\fivegoth
  \def\goth{\fam\gothfam\tengoth}
  
  \textfont\itfam=\tenit
  \def\it{\fam\itfam\tenit}
  \textfont\slfam=\tensl
  \def\sl{\fam\slfam\tensl}
  \textfont\bffam=\tenbf \scriptfont\bffam=\sevenbf
  \scriptscriptfont\bffam=\fivebf
  \def\bf{\fam\bffam\tenbf}
  \textfont\ttfam=\tentt
  \def\tt{\fam\ttfam\tentt}
  \abovedisplayskip=12pt plus 3pt minus 9pt
  \belowdisplayskip=\abovedisplayskip
  \abovedisplayshortskip=0pt plus 3pt
  \belowdisplayshortskip=4pt plus 3pt 
  \smallskipamount=3pt plus 1pt minus 1pt
  \medskipamount=6pt plus 2pt minus 2pt
  \bigskipamount=12pt plus 4pt minus 4pt
  \normalbaselineskip=12pt
  \setbox\strutbox=\hbox{\vrule height8.5pt depth3.5pt width0pt}
  \let\bigf@nt=\tenrm       \let\smallf@nt=\sevenrm
  \normalbaselines\rm}

\def\eightpoint{
  \textfont0=\eightrm \scriptfont0=\sixrm \scriptscriptfont0=\fiverm
  \def\rm{\fam\z@\eightrm}
  \textfont1=\eighti  \scriptfont1=\sixi  \scriptscriptfont1=\fivei
  \def\oldstyle{\fam\@ne\eighti}\let\old=\oldstyle
  \textfont2=\eightsy \scriptfont2=\sixsy \scriptscriptfont2=\fivesy
  \textfont\gothfam=\eightgoth \scriptfont\gothfam=\sixgoth
  \scriptscriptfont\gothfam=\fivegoth
  \def\goth{\fam\gothfam\eightgoth}
  
  \textfont\itfam=\eightit
  \def\it{\fam\itfam\eightit}
  \textfont\slfam=\eightsl
  \def\sl{\fam\slfam\eightsl}
  \textfont\bffam=\eightbf \scriptfont\bffam=\sixbf
  \scriptscriptfont\bffam=\fivebf
  \def\bf{\fam\bffam\eightbf}
  \textfont\ttfam=\eighttt
  \def\tt{\fam\ttfam\eighttt}
  \abovedisplayskip=9pt plus 3pt minus 9pt
  \belowdisplayskip=\abovedisplayskip
  \abovedisplayshortskip=0pt plus 3pt
  \belowdisplayshortskip=3pt plus 3pt 
  \smallskipamount=2pt plus 1pt minus 1pt
  \medskipamount=4pt plus 2pt minus 1pt
  \bigskipamount=9pt plus 3pt minus 3pt
  \normalbaselineskip=9pt
  \setbox\strutbox=\hbox{\vrule height7pt depth2pt width0pt}
  \let\bigf@nt=\eightrm     \let\smallf@nt=\sixrm
  \normalbaselines\rm}

\tenpoint

\def\pc#1{\bigf@nt#1\smallf@nt}         \def\pd#1 {{\pc#1} }

\catcode`\;=\active
\def;{\relax\ifhmode\ifdim\lastskip>\z@\unskip\fi
\kern\fontdimen2  -1.2 \fontdimen3 \string;}

\catcode`\:=\active
\def:{\relax\ifhmode\ifdim\lastskip>\z@\unskip\fi\penalty\@M\ \fi\string:}

\catcode`\!=\active
\def!{\relax\ifhmode\ifdim\lastskip>\z@
\unskip\fi\kern\fontdimen2  -1.1 \fontdimen3 \string!}

\catcode`\?=\active
\def?{\relax\ifhmode\ifdim\lastskip>\z@
\unskip\fi\kern\fontdimen2  -1.1 \fontdimen3 \string?}

\catcode`\«=\active 
\def«{\raise.4ex\hbox{%
 $\scriptscriptstyle\langle\!\langle$}}

\catcode`\»=\active 
\def»{\raise.4ex\hbox{%
 $\scriptscriptstyle\rangle\!\rangle$}}

\frenchspacing

\def\raggedbottom{\topskip 10pt plus 36pt\r@ggedbottomtrue}

\def\pointir{\unskip . --- \ignorespaces}

\def\Medbreak{\vskip-\lastskip\medbreak}

\long\def\th#1 #2\enonce#3\endth{
   \Medbreak\noindent
   {\pc#1} {#2\unskip}\pointir{\it #3}\smallskip}

\def\demonstration{\vskip-\lastskip\smallskip\noindent
 {\it Démonstration} : }

\def\decale#1{\smallbreak\hskip 28pt\llap{#1}\kern 5pt}
\def\decaledecale#1{\smallbreak\hskip 34pt\llap{#1}\kern 5pt}
\def\puce{\smallbreak\hskip 6pt{$\scriptstyle\bullet$}\kern 5pt}

\def\eqalign#1{\null\,\vcenter{\openup\jot\m@th\ialign{
\strut\hfil$\displaystyle{##}$&$\displaystyle{{}##}$\hfil
&&\quad\strut\hfil$\displaystyle{##}$&$\displaystyle{{}##}$\hfil
\crcr#1\crcr}}\,}

\catcode`\@=12

\showboxbreadth=-1  \showboxdepth=-1

\mathcode`A="7041 \mathcode`B="7042 \mathcode`C="7043 \mathcode`D="7044
\mathcode`E="7045 \mathcode`F="7046 \mathcode`G="7047 \mathcode`H="7048
\mathcode`I="7049 \mathcode`J="704A \mathcode`K="704B \mathcode`L="704C
\mathcode`M="704D \mathcode`N="704E \mathcode`O="704F \mathcode`P="7050
\mathcode`Q="7051 \mathcode`R="7052 \mathcode`S="7053 \mathcode`T="7054
\mathcode`U="7055 \mathcode`V="7056 \mathcode`W="7057 \mathcode`X="7058
\mathcode`Y="7059 \mathcode`Z="705A

\def\A{\mathord{\bf A}}
\def\A{\mathord{\bf A}}

\def\Card{\mathop{\rm Card}\nolimits}

\def\F{\mathord{\bf F}}
\def\F{\mathord{\bf F}}
\def\Gal{\mathop{\rm Gal}\nolimits}
\def\Gal{\mathop{\rm Gal}\nolimits}
\def\Hom{\mathop{\rm Hom}\nolimits}

\def\Kbar{\overline K}
\def\Kbar{{\overline K}}
\def\Ketoile{K^\times}
\def\Ketoile{K^{\times}}
\def\Kprimetoile{{\Kprim{}^{\times}}}
\def\kprimetoile{{\kprim{}^{\times}}}
\def\thetabarprim{\overline{\theta}'}
\def\Kprim{{K'}}
\def\Omegaprim{{\Omega'}}
\def\Mprim{{M'}}

\def\Lprim{{L'}}

\def\Pic{\mathop{\rm Pic}\nolimits}
\def\Pic{\mathop{\rm Pic}\nolimits}
\def\P{\mathord{\bf P}}
\def\P{\mathord{\bf P}}

\def\Xbar{\overline X}

\def\Xprim{{X'}}

\def\Z{\mathord{\bf Z}}
\def\Z{\mathord{\bf Z}}

\def\azeroo#1{A_0(#1)_0}
\def\azeroo#1{A_0(#1)_0}

\def\dii{\delta_2^2}
\def\dii{\delta_2^2}
\def\di{\delta_1^2}
\def\di{\delta_1^2}
\def\eibar{\underline{\ei}}

\def\eiibar{\underline{\eii}}

\def\eiiprim{\eii^\prime}
\def\eiiprim{e_2^\prime}
\def\eii{e_2}
\def\eii{e_2}
\def\eimeiibar{\underline{\ei-\eii}}
\def\eiprim{\ei^\prime}
\def\eiprim{e_1^\prime} 
\def\ei{e_1}
\def\ei{e_1} 

\def\epsii{\varepsilon_2} 
\def\epsii{\varepsilon_2} 
\def\epsi{\varepsilon_1}
\def\epsi{\varepsilon_1}
\def\eps{\varepsilon}
\def\eps{\varepsilon}

\def\ketoile{k^\times}
\def\kprim{{k'}}

\def\mgoth{{\goth m}}
\def\mgoth{{\goth m}}
\def\modm{(\mod\thinspace\mgoth)}
\def\modm{(\mod\thinspace\mgoth)}
\def\mod{\mathop{\rm mod.}\nolimits}
\def\mod{\mathop{\rm mod.}\nolimits}
\def\normlk{N_{L|K}(L^\times)}
\def\normlk{N_{L|K}(L^\times)}
\def\numero{n$^{\rm o}$}
\def\ocycle{\hbox{$0$-cycle}}
\def\ocycle{\hbox{$0$-cycle}}

\def\ogoth{{\goth o}}
\def\ogoth{{\goth o}}
\def\omegaprim{\omega^{\prime}}  

\def\qp{{\bf Q}_p}
\def\qp{{\bf Q}_p}

\def\thetabar{\overline{\theta}}

\def\underline#1{\omega(#1)}
\def\unites{\ogoth^\times}

\def\xbar{\underline{x}}

\def\xprim{x'}
\def\vprim{v'}
\def\hprim{h'}
\def\xprim{x^\prime}

\def\xxmeixmeiibar{\underline{x(x-\ei)(x-\eii)}}

\def\zmodiixzmodii{(\zmodii)^2}
\def\zmodiixzmodii{(\zmodii)^2}
\def\zmodii{\Z / 2\Z}
\def\zmodii{\Z / 2\Z}
\def\normlprimkprim{N_{L'|K'}}
\def\Nprim{N'}
\def\thetaprim{\theta'}
\def\Aprim{A'}
\def\dessus{\ufl{}{}{7mm}}

\def\hfl#1#2#3{\smash{\mathop{\hbox to#3{\rightarrowfill}}\limits
^{\textstyle#1}_{\textstyle#2}}}
\def\gfl#1#2#3{\smash{\mathop{\hbox to#3{\leftarrowfill}}\limits
^{\textstyle#1}_{\textstyle#2}}}

\def\qed{\raise -2pt\hbox{\vrule\vbox to 10pt{\hrule width 4pt
                 \vfill\hrule}\vrule}}
\def\cqfd{\unskip\penalty 500\quad\qed\medbreak}

\def\phi{\varphi}

\def\numero{n$^{\rm o}$}
\def\numeros{n$^{\rm os}$}

\def\droite#1{\,\hfl{#1}{}{8mm}\,}

\def\diagram#1{\def\normalbaselines{\baselineskip=0pt\lineskip=5pt}
\matrix{#1}}

\def\ufl#1#2#3{\llap{$\textstyle #1$}
\left\uparrow\vbox to#3{}\right.\rlap{$\textstyle #2$}}

\def\vide{\hbox{{\rm \O}}} 

\newcount\refno 

\long\def\ref#1:#2<#3>{                                        
\global\advance\refno by1\par\noindent                              
\llap{[{\bf\number\refno}]\ }{#1} \pointir{\it #2} #3\goodbreak }

\def\citer#1(#2){[{\bf\number#1}\if#2\empty\relax\else,\ #2\fi]}

\def\fleche{\rightarrow}

\newcount\numerodesection
\def\section#1{\bigbreak
 {\bf\number\numerodesection.\ \ #1}\nobreak\medskip
 \advance\numerodesection by1}

\newcount\numeroderemarque
\def\remarque{\advance\numeroderemarque by1\smallbreak
{\it Remarque\/}\ \number\numeroderemarque~:}

\newcount\formuleno
\def\numeroter{\global\advance\formuleno by1
 \leqno{(\number\formuleno)}}
\def\formule(#1){{\rm (\number#1)}}

\def\rem#1\endrem{%
\Medbreak
{\it#1\unskip} : }
\def\machin{\mathord{*}}
\def\classe#1{[\,#1\,]}

\def\somme{\mathop{\smash{\raise 2pt\hbox{$\sum$}}}\limits}

\def\long{\mathop{\rm long}\nolimits}

\newbox\bibbox
\setbox\bibbox\vbox{\bigbreak
\centerline{{\pc R{\'E}F{\'E}RENCES} {\pc BIBLIOGRAPHIQUES}}

\ref{\pc BLOCH} (S.):
Lectures on algebraic cycles,
<Durham~: Duke University Mathematics Department, 1980>
\newcount\blochlivre  \global\blochlivre=\refno

\ref{\pc CH{\^A}TELET} (F.):
Points rationnels sur certaines courbes et surfaces cubiques,
<Enseignement math., {\bf 5}, 1959, p.~153--170>
\newcount\chatelet  \global\chatelet=\refno

\ref{\pc COLLIOT}-{\pc TH{\'E}L{\`E}NE} (J-L.) et {\pc CORAY} (D. F.):
L'{\'e}quivalence rationnelle sur les points ferm{\'e}s des surfaces
rationnelles fibr{\'e}es en coniques,
<Compositio Math., {\bf 39} (3), 1979, p.~301--332>
\newcount\ctcoray  \global\ctcoray=\refno

\ref{\pc COLLIOT}-{\pc TH{\'E}L{\`E}NE} (J-L.) et {\pc SANSUC} (J-J.):
On the {C}how groups of certain rational surfaces: a sequel to
a paper of {S}.\ {B}loch,
<Duke Math.\ J., {\bf 48} (2), 1981, p.~421--447>
\newcount\ctsansucsequel  \global\ctsansucsequel=\refno

\ref{\pc COLLIOT}-{\pc TH{\'E}L{\`E}NE} (J-L.) et {\pc SANSUC} (J-J.):
La descente sur les vari{\'e}t{\'e}s rationnelles,
<dans Journ{\'e}es de G{\'e}om{\'e}trie alg{\'e}brique d'Anger, Alpen
aan den Rijn~: Sijthoff \& Noordhoff, 1980.>
\newcount\ctsansucdescente  \global\ctsansucdescente=\refno

\ref{\pc COOMBES} (K. R.) et {\pc MUDER} (D. J.):
Zero-cycles on del Pezzo surfaces over local fields,
<Journal of Algebra, {\bf 97}, 1985, p.~438--460.>
\newcount\coombesmuder  \global\coombesmuder=\refno

\ref{\pc CORAY} (D. F.) et {\pc TSFASMAN} (M. A.):
Arithmetic on singular del {P}ezzo surfaces,
<Proc.\ London Math.\ Soc.\ (3) {\bf 57} (1), 1988, p.~25--87.>
\newcount\coraytsfasman  \global\coraytsfasman=\refno

\ref{\pc DALAWAT} (C. S.):
Groupe des classes de $0$-cycles sur les surfaces
rationnelles d{\'e}finies sur un corps local,
<Th{\`e}se, Universit{\'e} de Paris-Sud, Orsay, 1993.>
\newcount\these  \global\these=\refno

\ref{\pc FULTON} (W.):
Intersection theory,
<Berlin~: Springer-Verlag, 1984.>
\newcount\fulton  \global\fulton=\refno

\ref{\pc MANIN} (Yu. I.):
Cubic forms,
<Amsterdam~: North-Holland Publishing Co., 1986.>
\newcount\manin  \global\manin=\refno

\ref{\pc SANSUC} (J-J.):
{\`A} propos d'une conjecture arithm{\'e}tique sur le groupe de
{C}how d'une surface rationnelle,
<S{\'e}minaire de Th{\'e}orie des nombres de Bordeaux, Expos{\'e} 33, 1982.>
\newcount\sansuc  \global\sansuc=\refno

} %\bibbox

\centerline{{\bf Le groupe de Chow d'une surface de Châtelet}}
\smallskip
\centerline{{\bf sur un corps local}}
\vskip\baselineskip
\centerline{Chandan Singh Dalawat}

\vskip 20pt plus10pt
\section{L'{\'e}nonc{\'e} des r{\'e}sultats}

Soit $K$ une extension de degr{\'e} fini du corps~$\qp$ des
nombres~\hbox{$p$-adiques} ($p$~nombre~premier~{\it impair\/}).
D{\'e}signons par~$\ogoth$ l'anneau des entiers de~$K$, par~$\mgoth$
l'ideal maximal de~$\ogoth$, par~$k=\ogoth/\mgoth$ le corps
r{\'e}siduel de~$K$ et par~$v:\Ketoile\rightarrow\Z$ la valuation
normalis{\'e}e de~$K$.  Soient~$d\in\Ketoile$ non~carr{\'e}
et~$\ei,\eii\in K$ deux {\'e}l{\'e}ments distincts non nuls.

Dans cette Note, on s'int{\'e}resse aux $K$-surfaces~$X$ qui sont lisses,
projectives et \hbox{$K$-birationnelles} {\`a} la surface affine
d'{\'e}quation
$$
y^2-dz^2=x(x-\ei)(x-\eii)\ ; \leqno{(1)}
$$ 
une telle surface est $K(\sqrt{d}\,)$-birationnelle au plan
projectif.  Dans la suite, on note $L=K(\sqrt{d}\,)$.

Les surfaces de Ch{\^a}telet \citer\chatelet(), \citer\manin() en
fournissent des exemples ; elles sont d{\'e}finies dans $\P({\cal
O}(2)\oplus{\cal O}(2)\oplus{\cal O}))$ (coordonn{\'e}es $y:z:t$)
au-dessus de la droite projective $\P_{1,K}$ (coordonn{\'e}es
$x:\xprim$) par l'{\'e}quation
$$
y^2-dz^2=x\xprim(x-\ei\xprim)(x-\eii\xprim)t^2.\leqno{(2)}
$$ 
Ce sont des fibr{\'e}s en coniques au-dessus de $\P_{1,K}$ avec quatre
fibres d{\'e}g{\'e}n{\'e}r{\'e}es au-dessus de $0,\ei,\eii$ et
$\infty$.  Ces $K$-mod{\`e}les lisses projectifs des
surfaces~(1) ont {\'e}t{\'e} construits dans \citer\ctsansucdescente().

{\`A} la permutation de $\ei,\eii$ pr{\`e}s, on peut supposer que
$v(\eii)\geq v(\ei)$. En outre, si l'on a $v(\eii)>v(\ei)$, le
changement de variable $x_1=x-\ei$  transforme
l'{\'e}quation~(1) en 
$$
y^2-dz^2=x_1(x_1-\eiprim)(x_1-\eiiprim)
\leqno{(3)}
$$ 
avec $\eiprim=-\ei$, \ $\eiiprim=\eii-\ei$ et l'on aura
$v(\eiiprim)=v(\eiprim)$ ; on peut donc supposer que $v(\eii)=v(\ei)$,
ce que nous faisons d{\'e}sormais.  On pose alors $r=v(\ei)=v(\eii)$.

Nous d{\'e}signons par $\azeroo{X}$ le groupe des \ocycle s de
degr{\'e} $0$ sur $X$, modulo l'{\'e}quivalence rationnelle
\citer\blochlivre(), \citer\fulton(), \citer\manin().  Divers auteurs
\citer\ctsansucsequel(), \citer\sansuc() ont cherch{\'e}
{\`a} calculer explicitement ce groupe.  Un r{\'e}sultat d{\^u} {\`a}
Colliot-Thélène \citer\coraytsfasman(prop.~4.7) se r{\'e}sume ainsi~:

\th PROPOSITION 1 (Colliot-Thélène)
\enonce
Supposons que l'extension $L$ est non
ramifi{\'e}e. Pour toute $K$-surface $X$ lisse projective
$K$-birationnelle {\`a} la surface~$(1)$, le groupe $\azeroo{X}$ est
isomorphe
$$\vbox{\halign{\hfil\it #\/\rm)\ &{\`a}#&\qquad\hss$#$\hss\qquad
&{\it si\/ #}\hfil\cr
  i&&0^{\phantom{2}} & $r$ est pair et\/ $v(\ei-\eii)=r$,\cr
 ii&&\zmodii^{\phantom{2}} & $r$ est pair et\/ $v(\ei-\eii)>r$,\cr
iii&&\zmodiixzmodii & $r$ est impair.\cr}}
$$
\endth

L'{\'e}nonc{\'e} de \citer\coraytsfasman(prop.~4.7) normalisait les
valuations de $\ei$ et $\eii$ par $v(\ei)=0\ \hbox{ou}\ 1$ et
$v(\eii)\ge0$, ce qui amenait {\`a} consid{\'e}rer sept cas de figure.
Une d{\'e}monstration simplifi{\'e}e est donn{\'e}e au \numero~4.  En
effet, nous ne distinguons que les trois cas de l'{\'e}nonc{\'e}
(visiblement irr{\'e}ductible de l'un {\`a} l'autre).

Notre propos ici est de d{\'e}terminer le groupe $\azeroo{X}$ lorsque
l'extension $L$ est ramifi{\'e}e (Proposition~2) et de calculer
l'homomorphisme de restriction $\azeroo{X}\fleche\azeroo{\Xprim}$,
o{\`u} $\Kprim$ est une extension de degr{\'e} fini de~$K$ et
$\Xprim=X\times_K\Kprim$ (Proposition~3).

Les deux principaux r{\'e}sultats d{\'e}montr{\'e}s ici sont donc les
suivants~:\nobreak

\th PROPOSITION 2
\enonce
Supposons que l'extension\/ $L$ est ramifi{\'e}e.  Pour
toute\/ $K$-surface\/ $X$ lisse projective\/ $K$-birationnelle {\`a}
la surface\/~$(1)$, le groupe\/ $\azeroo{X}$ est isomorphe
$$
\vbox{\halign{\hfil\it #\/\rm)\ &{\`a}#&\qquad\hss$#$\hss\qquad&
{\it si\/ #}\hfil\cr
  i&&\zmodii^{\phantom{2}} & $\ei/\eii\equiv1\,\modm$ et\/ $\ei\in\normlk$,\cr 
 ii&&\zmodiixzmodii & $\ei/\eii\equiv1\,\modm$ et\/ $\ei\notin\normlk$,\cr 
iii&&\zmodiixzmodii & $\ei/\eii\not\equiv1\,\modm$.\cr}}
$$
\endth

La d{\'e}monstration est donn{\'e}e aux \numeros~5--8.  (Noter que les
trois conditions sont en fait sym{\'e}triques en $\ei$ et $\eii$.  Par
ailleurs, les types~{\it ii\/}) et~{\it iii\/}) sont
$K$-birationnellement distincts (voir {\it Remarque\/}~10).)

\remarque\ Soient $a$ et $c$ deux {\'e}l{\'e}ments
de $K$, distincts de $0$ et de $1$.  On sait
\citer\coraytsfasman(prop.~4.6) que la surface $X$ d{\'e}finie par le
syst{\`e}me d'{\'e}quations 
$$
\cases{a(d\xi^2-\zeta^2)=(\tau+\eta)(\eta+\nu),\cr
c(d\xi^2-\eta^2)=(\tau+\zeta)(\zeta-\nu),\cr}\leqno{(4)}
$$
(une surface de del Pezzo de degr{\'e}~4) est $K$-birationnelle {\`a}
la surface~(1) avec $\ei=ac$ et $\eii=a+c-1$.  Par l'invariance
birationnelle du groupe de Chow de $0$-cycles de degr{\'e}\/~$0$
\citer\ctcoray(prop.~6.3) (cf.\ \numero~3, th.~1), les propositions~1
et~2 d{\'e}terminent le groupe $\azeroo{X}$.  Le r{\'e}sultat
principal de Coombes et Muder
\citer\coombesmuder(th{\'e}or{\`e}mes~4.4 et~4.5) calcule ce groupe
dans le cas o{\`u} $L$ est non ramifi{\'e}e, par une {\'e}tude du
$\Gal(\Kbar\,|\,K)$-ensemble des droites trac{\'e}es sur~(4), $\Kbar$
{\'e}tant une cl{\^o}ture alg{\'e}brique de~$K$ (voir aussi
\citer\these(exemple~6.3)).

\goodbreak

\th PROPOSITION 3
\enonce 
Soient\/ $\Kprim$ une extension finie de~$K$, $X$ une $K$-surface
lisse projective $K$-birationnelle {\`a} la surface~$(1)$, et\/
$\Xprim=X\times_K\Kprim$.  L'homomorphisme de restriction
$\azeroo{X}\fleche\azeroo{\Xprim}$ est trivial si le degr{\'e}\/
$n=[\Kprim:K]$ est pair ; c'est un isomorphisme sinon.
\endth

La d{\'e}monstration est donn{\'e}e au \numero~9.

\remarque\ L'homomorphisme de corestriction
$\azeroo{\Xprim}\fleche\azeroo{X}$ est toujours injectif pour de
telles surfaces, quel que soit le degr{\'e}~$n$.  C'est m{\^e}me un
isomorphisme pour $n$~impair, d'apr{\`e}s les propositions~1--3 et le
fait que le compos{\'e}
$\azeroo{X}\fleche\azeroo{\Xprim}\fleche\azeroo{X}$ est la
multiplication par~$n$ (voir par exemple \citer\fulton()).  Si $n$ est
pair et si $\Kprim$ contient $\sqrt{d}$, le groupe $\azeroo{\Xprim}$
est nul car alors $\Xprim$ est $\Kprim$-birationnelle {\`a} $\P_{2}$
\citer\blochlivre(p.~7.1).  L'{\'e}nonc{\'e} g{\'e}n{\'e}ral
r{\'e}sulte d'une proposition communiqu{\'e}e {\`a} l'auteur par
Colliot-Thélène ; elle permet par ailleurs de retrouver la prop.~3.

\goodbreak

{\it Je remercie Jean-Louis Colliot-Thélène pour m'avoir initi{\'e} {\`a} ce
sujet et pour ses conseils.  Lorsque ce travail a {\'e}t{\'e}
commenc{\'e}, l'auteur b{\'e}n{\'e}ficiait d'une bourse du
Gouvernement fran{\c c}ais.}

\goodbreak

\section{Des formulations {\'e}quivalentes}

Une reformulation des propositions~1 et~2 fait voir plus clairement
comment le groupe $\azeroo{X}$ d{\'e}pend du type de r{\'e}duction de
la surface~(1). 

Supposons que $L$ est non ramifi{\'e}e.  Soient $\pi$ une
uniformisante de $K$ et $\eps_i$ les images dans $k$ des
$e_i/\pi^{2m}$ si $r=2m$ ou $r=2m+1$ (cf.\ {\it Remarque\/}~5).

Consid{\'e}rons la $k$-cubique plane $C$ et le point $P$ sur $C$~:
$$
C : y^2=x(x-\epsi)(x-\epsii)\ ;\ \
P=(\epsi,0) .\leqno{(5)} 
$$

La courbe $C$ ne d{\'e}pend pas du choix de $\pi$, {\`a}
$k$-isomorphisme pr{\`e}s.  Elle est lisse si $\epsi\neq\epsii$ et
singuli{\`e}re en $P$ sinon ; $P$ est un point double ordinaire si
$\epsi=\epsii\neq0$ et un point de rebroussement si $\epsi=\epsii=0$.

\th PROPOSITION 4 
\enonce 
Supposons que l'extension\/ $L$ est non ramifi{\'e}e.  Pour toute
$K$-surface $X$ lisse projective $K$-birationnelle {\`a} la
surface~$(1)$, le groupe $\azeroo{X}$ est isomorphe
$$\vbox{\halign{\hfil\it #\/\rm)\ &{\`a}#&\qquad\hss$#$\hss\qquad&
{\it si\/ #}\hfil\cr
  i&&0^{\phantom{2}} & $P$ est un point r{\'e}gulier,\cr
 ii&&\zmodii^{\phantom{2}} & $P$ est un point double ordinaire,\cr 
iii&&\zmodiixzmodii & $P$ est un point de rebroussement.\cr}}
$$
\endth
Il est clair que cette proposition est {\'e}quivalente {\`a} la prop.~1.

\goodbreak

Supposons maintenant que l'extension\/ $L$ est ramifi{\'e}e.  Soit
$\pi$ une uniformisante de $K$ telle que $\pi\in\normlk$ et
d{\'e}signons par $\omega:\Ketoile\fleche\ketoile$ l'homomorphisme
$x\longmapsto x/{\pi^{v(x)}}\,\modm$.  Les normes $\normlk$ sont
pr{\'e}cis{\'e}ment les $x\in \Ketoile$ pour lesquels
$\xbar\in\ketoile$ est un~carr{\'e} (voir \numero~5).

Soient $\eps_i=\omega(e_i)$.  La courbe $C$~(5) est ind{\'e}pendante du
choix de $\pi$, {\`a} $k$-isomorphisme pr{\`e}s, puisque
$v(\ei)=v(\eii)$.  Elle est lisse si $\epsi\neq\epsii$ et
singuli{\`e}re au point $P$ sinon ; c'est alors un point double
ordinaire qui est d{\'e}ploy{\'e} si et seulement si l'{\'e}l{\'e}ment
$\epsi=\epsii$ de $\ketoile$ est un carr{\'e}.

Une fa{\c c}on succincte d'{\'e}noncer la prop.~2 est la suivante~:

\th PROPOSITION 5
\enonce
Supposons que l'extension\/ $L$ est ramifi{\'e}e.  Pour toute\/
$K$-surface $X$ lisse projective\/ $K$-birationnelle {\`a} la
surface\/~$(1)$, le groupe\/ $\azeroo{X}$ est isomorphe {\`a}\/
$\zmodii$ si le point\/ $P$ de la cubique\/~$C$~$(5)$ est un point
double ordinaire d{\'e}ploy{\'e} ; il est isomorphe {\`a}\/
$\zmodiixzmodii$ sinon.
\endth
\rem D{\'e}monstration de l'{\'e}quivalence avec la prop.~{\rm 2}\endrem
Dans les cas~{\it i\/}) et {\it ii\/}), cette courbe est
singuli{\`e}re au point~$(\epsi,0)$ ; avec le changement de variable
$x_1=x-\epsi$ pour que la singularit{\'e} soit {\`a}~l'origine, le
terme homog{\`e}ne de plus bas~degr{\'e} du polyn{\^o}me
d{\'e}finissant~$C$ est {\'e}gal~{\`a} $y^2-\epsi^{\phantom{2}}
x_1^2$~; c'est donc un produit de deux facteurs lin{\'e}aires sur~$k$
si et seulement si $\epsi$ est un~carr{\'e} dans~$\ketoile$, ce qui
{\'e}quivaut {\`a} dire que~$\ei\in\normlk$.  Enfin, dans le cas~{\it
iii\/}), la cubique~$C$ est lisse.\cqfd

\section{La m{\'e}thode de calcul}

Elle est calqu{\'e}e sur la m{\'e}thode de
\citer\coraytsfasman(prop.~4.7) et repose sur les quatre
th{\'e}or{\`e}mes suivants~:

\th
TH{\'E}OR{\`E}ME 1 (Colliot-Thélène et Coray \citer\ctcoray(prop.~6.3)) 
\enonce
Le groupe de Chow des $0$-cycles de degr{\'e}\/~$0$ sur une
surface\/~$S$ lisse, projective, absolument connexe sur un
corps\/~$F$ est un invariant\/ $F$-birationnel de\/~$S$.
\endth

Ceci nous permet de ne d{\'e}montrer les propositions~1 et~2 que pour
les surfaces de Ch{\^a}telet~; on suppose d{\'e}sormais que $X$
d{\'e}signe une telle surface, donn{\'e}e par l'{\'e}quation~(2).

Notons $O$ le point singulier de la fibre {\`a} l'infini de la
fibration en coniques $f:X\fleche\P_1$ et consid{\'e}rons l'application
$$\gamma:X(K)\fleche\azeroo{X},\quad \gamma(Q)=Q-O.\leqno{(6)}$$
\th
TH{\'E}OR{\`E}ME 2 (Colliot-Thélène et Coray \citer\ctcoray(th.~C)) 
\enonce
Pour les surfaces de Ch{\^a}telet sur un corps local num{\'e}rique,
l'application~$\gamma$ est surjective.
\endth
\goodbreak

Soit $\Kbar$ une cl{\^o}ture alg{\'e}brique de $K$ et posons
$\Xbar=X\times_{K}\Kbar$.  Le module galoisien $\Pic(\Xbar)$ est un
groupe commutatif libre de type fini.  Notons $S$ le $K$-tore dont le
$\Gal(\Kbar\,|\,K)$-module de caract{\`e}res est le groupe
$\Pic(\Xbar)$ ; on a $S(\Kbar)=\Hom_{\Z}(\Pic(\Xbar),\,\Kbar^\times)$.
Colliot-Thélène et Sansuc ont construit l'homomorphisme caract{\'e}ristique
$$\phi:\azeroo{X}\fleche H^1(K,S(\Kbar))\leqno{(7)}$$
pour la d{\'e}finition duquel le lecteur est pri{\'e} de se reporter
{\`a} \citer\ctsansucsequel().

\th
TH{\'E}OR{\`E}ME 3 (Colliot-Thélène et Sansuc \citer\ctsansucsequel(th.~5)) 
\enonce
Pour les surfaces fibr{\'e}es en coniques au-dessus de la droite
projective, l'homomorphisme~$\phi$ est injectif.
\endth

Une {\'e}tude du $\Gal(\Kbar\,|\,K)$-module $\Pic(\Xbar)$ fournit un
isomorphisme (cf.~\citer\sansuc())
$$\iota:H^1(K,S(\Kbar))\fleche(\Ketoile\!/\,\normlk)^2.\leqno{(8)}$$
Notons $h:\Ketoile\fleche\zmodii$ l'homomorphisme surjectif de noyau
$\normlk$ ; on identifie ainsi le quotient $\Ketoile\!/\,\normlk$
{\`a} $\zmodii$.

La th{\'e}orie de la descente de \citer\ctsansucdescente()
(cf.~\citer\sansuc(\numero~III, p.~33-11) et
\citer\coraytsfasman(p.~59)) donne une description explicite de la
compos{\'e}e des applications (6)--(8) :
\th
TH{\'E}OR{\`E}ME 4 (Colliot-Thélène et Sansuc \citer\ctsansucdescente({\numero~IV})) 
\enonce
L'application compos{\'e}e $\theta=\iota\circ\phi\circ\gamma:
X(K)\fleche\zmodiixzmodii$ est donn{\'e}e par la formule~: 
$$\def\classe#1{h(#1)}
\theta(y:z:t\;;\,x) =\cases{
(\classe{1},\, \classe{1}) =0 & si $x=\infty$,\cr
(\classe{\ei\eii},\, \classe{-\ei}) & si $x=0$,\cr
(\classe{\ei},\, \classe{\ei(\ei-\eii})& si $x=\ei$,\cr
(\classe{x},\, \classe{x-\ei}) & sinon.\cr}\leqno{(9)}
$$
\endth

Comme tous les points d'une fibre de $f : X(K)\fleche\P_1(K)$ sont des
\ocycle s rationnellement {\'e}quivalents, on a une application
induite $f(X(K))\fleche\zmodiixzmodii$ encore not{\'e}e~$\theta$.
D'apr{\`e}s les th{\'e}or{\`e}mes~2 et~3, l'homomorphisme
$\iota\circ\phi$ identifie le groupe $\azeroo{X}$ {\`a}
$\theta(f(X(K)))$.

La partie $f(X(K))$ de $\P_1(K)$ contient exactement
$0,\,\ei,\,\eii,\, \infty$ et les $x\in\Ketoile$
satisfaisant $x(x-\ei)(x-\eii)\in\normlk$ (cf.\ {\'e}quation~(1)).

\remarque\ Soit $N\subset f(X(K))$ la partie contenant $0$ et
les $x\in f(X(K))$ tels que $v(x)=r$. Pour calculer l'image de
$f(X(K))$ par $\theta$, on peut se restreindre {\`a} $N$ ({\it i.e.}
$\theta(N)$ engendre le groupe $\iota\circ\phi(\azeroo{X})$) gr{\^a}ce
au lemme suivant :

\th
LEMME 1
\enonce
On a $\theta(x)=0$ si $v(x)<r$ et\/ $\theta(x)=\theta(0)$ si $v(x)>r$.
\endth
\demonstration        
Comme, pour $x\in f(X(K))$ (distinct de $0,\,\ei,\,\eii$ et
$\infty$),\break on a $h(x)+h(x-\ei)+h(x-\eii)=0$, il suffit de
v{\'e}rifier que\break $h(x-e_i)=h(x)$ si $v(x)<r$ et
$h(x-e_i)=h(-e_i)$ si $v(x)>r$.\break C'est clair si l'extension $L$
est non ramifi{\'e}e car $h$ est alors le compos{\'e}
$\def\droite#1{\hfl{#1}{}{7mm}}\Ketoile\droite{v}\Z\fleche\zmodii$ ;
lorsque $L$ est ramifi{\'e}e, ceci r{\'e}sulte de ce que $h$ se
factorise en
$\Ketoile\fleche\unites\fleche\ketoile\fleche\ketoile\!/k^{\times2}$,
le premier homomorphisme {\'e}tant $x\longmapsto
x/{\pi^{v(x)}}\,\modm$, pour une uniformisante $\pi$ de $K$ telle que
$\pi\in\normlk$ (cf.\ \numero~5).\cqfd

En derni{\`e}re analyse, d{\'e}montrer les propositions~1 et~2
revient {\`a} d{\'e}terminer $\theta(N)$.

\goodbreak

\section{Le cas non ramifi{\'e}}

Dans ce \numero, on d{\'e}montre la prop.~1.  Supposons donc que
l'extension\/ $L$ est non ramifi{\'e}e.

La partie $N$ de $K$ ({\it Remarque\/}~3) contient pr{\'e}cis{\'e}ment
$0,\ei,\eii$ et tous les $x\in K$ (distincts de $0,\,\ei$ et $\eii$)
satisfaisant $v(x)=r$ pour lesquels $v(x(x-\ei)(x-\eii))$ est pair.
Il s'agit de calculer $\theta(N)$.

\rem Le cas i) $r$ est pair et $v(\ei-\eii)=r$\endrem
L'image de $0,\,\ei$ et $\eii$ par $\theta$ est triviale ; pour tout
autre $x\in N$, on a $\theta(x)=0$ si $v(x-\ei)=r$ ; supposons donc
que $n=v(x-\ei)$ est $>r$.  Comme $x-\eii=x-\ei+(\ei-\eii)$, on a
$v(x-\eii)=r$ et comme $x\in N$, \ $n$ doit {\^e}tre pair, d'o{\`u}
$\theta(x)=0$.

\rem Le cas ii) $r$ est pair et $v(\ei-\eii)>r$\endrem On pose
$n=v(\ei-\eii)$.  Si $n$ est impair, on a $\theta(\ei)=(0,\,1)$ et la
preuve est termin{\'e}e car l'image de tout {\'e}l{\'e}ment de $N$ est
de la forme $(0,\machin)$.  Supposons donc que $n$ est pair et soit
$m$ un entier impair satisfaisant $r<m<n$ ; on pose $x=\ei+\pi^m$,
o{\`u} $\pi$ est une uniformisante de $K$. Or $x$ appartient {\`a} $N$
puisque $v(x)=r$, $v(x-\ei)=m$ et $v(x-\eii)=m$ ; on a
$\theta(x)=(0,1)$.

\rem Le cas iii) $r$ est impair\endrem
On a $\theta(0)=(0,1)$ et $\theta(\ei)=(1,\machin)$. \cqfd

(Il est {\`a} noter que l'hypoth{\`e}se ``~$p\neq2$~'' n'a pas
{\'e}t{\'e} utilis{\'e}e.)

\section{Uniformisantes adapt{\'e}es}

Les \numeros~5--8 sont consacr{\'e}s {\`a} la d{\'e}monstration de la
prop.~2.

Supposons donc que l'extension $L$ est ramifi{\'e}e.  Dans le
pr{\'e}sent \numero, on attache, au choix~$\pi$ d'une uniformisante
de~$K$ telle que $\pi\in\normlk$, une application de r{\'e}duction
$\omega:K\fleche k$.  

\goodbreak
Soit $M$ la partie de $N$ ({\it Remarque\/}~3)
contenant pr{\'e}cis{\'e}ment $\ei$ et les $x\in N$ tels que
$\xbar\neq\eibar$.  

Dans le \numero~6, on caract{\'e}rise l'image
$\omega(M)\subset k$.  Ensuite, au \numero~7, on d{\'e}finit une
certaine application $\thetabar: k\fleche\zmodiixzmodii$~(11) et l'on
d{\'e}termine le sous-groupe engendr{\'e} par~$\thetabar(\omega(M))$.
Bien que les applications $\theta$ et $\thetabar\circ\omega$ ne
commutent pas sur $N$, elles le font sur $M$.  Par bonheur, on
constate que $\theta(N)=\theta(M)$~(\numero~8), achevant ainsi de
d{\'e}montrer la prop.~2.

\th D{\'E}FINITON 1
\enonce
On dira qu'une uniformisante $\pi$ de $K$ est adapt{\'e}e {\`a\/} $L$
si $\pi\in\normlk$.
\endth

\th LEMME 2
\enonce
Il existe toujours des uniformisantes adapt{\'e}es {\`a}\/ $L$.  Pour
une telle uniformisante $\pi$, les normes $\normlk$ sont
pr{\'e}cis{\'e}ment les $x\in \Ketoile$ pour lesquels
$x/{\pi^{v(x)}}\,\modm$ est un~carr{\'e} dans $\ketoile$.
\endth
\goodbreak
\demonstration 
Soient $\varpi$ une uniformisante de $K$ et $t$ un g{\'e}n{\'e}rateur
du groupe $\ketoile\!/k^{\times2}$.  {\`A} l'aide de $\varpi$, on
identifie $\Ketoile\!/\,K^{\times2}$ {\`a}
$\zmodii\times\ketoile\!/k^{\times2}$.  Alors $\normlk/K^{\times2}$
s'identifie soit {\`a} $\zmodii\times\{1\}$, soit {\`a} $\{(0,1),\
(1,t)\}$; dans le premier cas, l'uniformisante $\varpi$ est
adapt{\'e}e {\`a} $L$ et dans le second $u\varpi$ ($u\in\unites$
relevant $t$) l'est.\cqfd

Par exemple, lorsque $K=\qp$ et $d=p$, on peut prendre
$\pi=p$ si $p\equiv1\,(\mod 4)$ et $\pi=-p$ si $p\equiv-1\,(\mod 4)$.

Choisissons d{\'e}sormais une uniformisante $\pi$ de~$K$ adapt{\'e}e
{\`a}~$L$ et d{\'e}signons par $\omega:\Ketoile\fleche\ketoile$
l'homomorphisme $x\mapsto x/{\pi^{v(x)}}\,\modm$.  En posant
$\omega(0)=0$, prolongeons-le en une application $\omega :K\fleche k$.

\remarque\ L'homomorphisme compos{\'e}
$\Ketoile\fleche\ketoile\fleche\ketoile\!/k^{\times2}$ est
ind{\'e}pendant du choix de~$\pi$ car son noyau est {\'e}gal
{\`a}~$\normlk$ (lemme~2) ; c'est l{\`a} l'homomorphisme~$h$ qui
intervient dans~(9).

\remarque\ Dans l'{\'e}quation~(1), on peut supposer sans perte de
g{\'e}n{\'e}ralit{\'e} que $d,\,\ei$ et $\eii$ appartiennent {\`a}
$\ogoth$, que $v(d)=0$~ou~$1$ et que $r=0$~ou~$1$.  Si $v(d)=1$, on
peut m{\^e}me supposer que $r=0$ car l'extension $L=K(\sqrt{d}\,)$ est
alors ramifi{\'e}e.  {\'E}crivons en effet $e_i=\pi^r u_i$, avec
$u_i\in\unites$ et $\pi$ une uniformisante de $K$ adapt{\'e}e {\`a}
$L$.  Soient $a,b\in\Ketoile$ tels que $\pi^{-3r}=a^2-db^2$ (ce qui
est possible puisque $\pi\in\normlk\,$) ; le changement de variables
$y_1=ay+dbz$, \ $z_1=az+by$, \ $x_1=\pi^{-r}x$ ram{\`e}ne (1) {\`a}
$$
y_1^2-dz_1^2=x_1(x_1-u_1)(x_1-u_2).\leqno{(10)}
$$ 
Par contre, lorsque l'extension $L$ est non ramifi{\'e}e et $r=1$, il
n'existe pas de fonctions $x_1,y_1,z_1\in K(X)$ et
d'unit{\'e}s $u_1,u_2\in\unites$ satisfaisant la relation~(10) car le
groupe de Chow est diff{\'e}rent dans le cas~{\it iii\/}) de la
prop.~1 de ce qu'il est dans les deux autres cas.
\goodbreak
\th LEMME 3
\enonce
Si\/ $x_1, x_2\in\Ketoile$ ont la m{\^e}me valuation et si\/
$\omega(x_1)\neq\omega(x_2)$, on a\/ $\omega(x_1-x_2)=\omega(x_1)-\omega(x_2)$.
\endth
\demonstration Posons $n=v(x_i)$ et {\'e}crivons $x_i=u_i\pi^n$ ($u_1,
u_2\in\unites$) ; on a $\omega(x_i)=u_i\,\modm$.  Puisque
$\omega(x_1)\neq\omega(x_2)$, $u_1-u_2$ est une unit{\'e} et par suite
$v(x_1-x_2)=n$. D'o{\`u}~: $\omega(x_1-x_2) =u_1-u_2\,\modm$.\cqfd

\goodbreak

\section{R{\'e}duction {\`a} un probl{\`e}me combinatoire}

Rappelons que $M\subset N$ (\numero~5) contient pr{\'e}cis{\'e}ment
$\ei$ et les $x\in N$ tels que $\xbar\neq\eibar$ ; on a clairement
$\omega(M)=\omega(N)$.

On pose $\eps_i=\omega(e_i)$ ; on a
$\epsii\in\omega(M)$, m{\^e}me si $\eii\notin M$.
\goodbreak
\th LEMME 4
\enonce
La partie\/ $\omega(M)\subset k$ contient pr{\'e}cis{\'e}ment\/ tous les
$t\in k$ pour lesquels\/ $t(t-\epsi)(t-\epsii)$ est
un~carr{\'e}.
\endth
(Autrement dit, $\omega(M)$ est l'image de l'ensemble des points
rationnels $C(k)$ sur la cubique~(5) par le morphisme
$x:C\fleche\A^1_k$.)
\demonstration
Si $\xbar$ (pour $x\in M$) est distinct de $0,\,\epsi$ et
$\epsii$, on~a $\omega(x-e_i)=\omega(x)-\omega(e_i)$ (lemme~3).  Comme
le premier membre de
$$
\xxmeixmeiibar= \xbar(\xbar-\epsi)(\xbar-\epsii)
$$
est un~carr{\'e}, \ $\xbar$ v{\'e}rifie la condition de l'{\'e}nonc{\'e}.

Inversement, si $t\in k$ (distinct de $0,\,\epsi$ et $\epsii$) v{\'e}rifie
cette conditon, on pose $x=u\pi^{r}$ avec
$u\in\unites$ relevant $t$ ; on a
$\omega(x-e_i)=t-\eps_i$ par le lemme~3, si bien que
$x(x-\ei)(x-\eii)\in\normlk$, c'est-{\`a}-dire que $x\in M$ ;
on a $\xbar=t$.\cqfd

\th LEMME 5  
\enonce       
Si\/ $\epsi=\epsii$, \ $\omega(M)$ contient tous les carr{\'e}s de\/
$\ketoile$, et aucun non carr{\'e}, {\`a} l'exception {\'e}ventuelle
de\/ $\epsi=\epsii$.
\endth        
\demonstration C'est clair, d'apr{\`e}s le crit{\`e}re du lemme~4.\cqfd

\section{R{\'e}solution du probl{\`e}me combinatoire}

Dans ce \numero, on fait abstraction de l'origine du probl{\`e}me ;
curieusement, cela le rend plus concret.

Soit $k$ un corps fini de caract{\'e}ristique impair.  On fixe, une
fois pour toutes, un {\'e}l{\'e}ment {\it non carr{\'e}\/}
$\alpha\in\ketoile$.

\th LEMME 6
\enonce
Pour\/ $\eps\in\ketoile$, il existe toujours un\/ $\xi\in\ketoile$ tel
que\/ $\xi^2-\eps$ soit non~carr{\'e}, sauf si\/ $k=\F_3$ et\/
$\eps=1$.
\endth
\demonstration
Consid{\'e}rons la conique $\xi^2-\eps\zeta^2 =\alpha\eta^2$ dans le
plan $\P_{2,k}$ de coordonn{\'e}es homog{\`e}nes $\xi : \eta :
\zeta$~; elle poss{\`e}de $q+1$~points $k$-rationnels (o{\`u}
$q=\Card(k)$) dont aucun n'est sur la droite~$\zeta=0$ et aucun sur la
droite~$\eta=0$ non plus si $\eps$ n'est pas un carr{\'e}.\cqfd

Les trois lemmes qui suivent sont des cons{\'e}quences directes de ce
que tout espace homog{\`e}ne principal $E$ sous une $k$-courbe
ab{\'e}lienne est trivial, c'est-{\`a}-dire que $E(k)\neq\vide$.  En
effet, dans chaque cas, les {\'e}quations {\`a} r{\'e}soudre (rendues
homog{\`e}nes) d{\'e}finissent une courbe~$E$ lisse, projective,
absolument connexe, intersection de deux quadriques dans $\P_3$ ;
c'est donc un espace homog{\`e}ne principal sous sa jacobienne.  Les
hypoth{\`e}ses suppl{\'e}mentaire sont l{\`a} pour s'assurer que les
points rationnels de cette courbe sont tous dans l'ouvert
$\xi\eta\zeta\lambda\neq0$ ; si $(\xi:\eta:\zeta:1)\in E(k)$, \
$\xi,\eta,\zeta$ r{\'e}pondent {\`a} la question.

\th LEMME 7
\enonce
Soient\/ $\di$ et\/ $\dii$ des {\'e}l{\'e}ments distincts de\/
$\ketoile$.  Il existe $\xi,\,\eta,\,\zeta\in\ketoile$ tels que\/
$\alpha\xi^2-\di=\eta^2$ et\/ $\alpha\xi^2-\dii=\alpha\zeta^2$.
\endth
\th LEMME 8
\enonce
Soient\/ $\di$ et\/ $\alpha\dii$ des {\'e}l{\'e}ments de\/
$\ketoile$. Supposons que\/ $\di-\alpha\dii$ est un carr{\'e}.  Il
existe alors\/ $\xi,\,\eta,\,\zeta\in\ketoile$ tels que\/
$\xi^2-\di=\alpha\eta^2$ et\/ $\xi^2-\alpha\dii=\alpha\zeta^2$.
\endth
\th LEMME 9
\enonce
Soient\/ $\alpha\di$ et\/ $\alpha\dii$ des {\'e}l{\'e}ments distincts
de\/ $\ketoile$.  Supposons que\/ $-1$ n'est pas un carr{\'e} dans\/
$\ketoile$.  Il existe alors\/ $\xi,\,\eta,\,\zeta\in\ketoile$ tels que\/
$\xi^2-\alpha\di=\alpha\eta^2$ et\/ $\xi^2-\alpha\dii=\alpha\zeta^2$.
\endth
Comme il a {\'e}t{\'e} dit, ces trois lemmes r{\'e}sultent de ce que
les {\'e}quations qui interviennent, rendues homog{\`e}nes,
d{\'e}finissent une courbe de genre~1 sur le corps fini~$k$ dont aucun
point $k$-rationnel n'est sur les axes de coordonn{\'e}es.

\bigbreak
Soient maintenant $\epsi, \epsii$ des {\'e}l{\'e}ments (distincts ou
non) de $\ketoile$ et $\eps_{1,2}=\epsi-\epsii$ si $\epsi\neq\epsii$,
arbitraire dans $\ketoile$ sinon.  D{\'e}finissons une application
$\thetabar :k\fleche\zmodiixzmodii$ par la formule ({\`a} comparer
avec (9)) :
$$
\thetabar(t) =\cases{ %(\classe{1},\, \classe{1})=0 & si $t=\infty$,\cr
(\classe{\epsi\epsii},\, \classe{-\epsi}) & si $t=0$,\cr
(\classe{\epsi},\, \classe{\epsi\eps_{1,2}})& si $t=\epsi$,\cr 
(\classe{t},\, \classe{t-\epsi}) & sinon,\cr}\leqno{(11)}
$$
o{\`u} $\classe{\phantom{x}} : \ketoile\fleche\zmodii$ est le quotient
modulo les carr{\'e}s.

Soit $\Omega$ la partie de $k$ contenant tous les $t\in k$ pour
lesquels $t(t-\epsi)(t-\epsii)$ est un~carr{\'e} (cf.~lemme~4).

\th LEMME 10
\enonce
Si\/ $\epsi=\epsii$ et si c'est un carr{\'e}, alors\/
$\thetabar(\Omega)$ engendre le sous-groupe $\{0\}\times\zmodii$.
\endth
\demonstration
L'ensemble $\Omega$ contient tous les carr{\'e}s de $\ketoile$ et ne
contient aucun non carr{\'e} (cf.~lemme~5) ; il est donc clair qu'un
{\'e}l{\'e}ment du type $(1, \machin)$ ne peut pas appartenir {\`a}
$\thetabar(\Omega)$.

Si $k=\F_3$ et $\epsi=1$, on a $\thetabar(0)=(0,1)$ (comme quoi, le
petit z{\'e}ro sert~{\`a} quelque chose).  Sinon, le lemme~6 fournit un
$\xi\in\ketoile$ tel que $\xi^2-\epsi$ soit non~carr{\'e} ; on a alors
$\xi^2\in\Omega$ et  $\thetabar(\xi^2)=(0,1)$.\cqfd
\goodbreak
\th LEMME 11
\enonce
Si\/ $\epsi=\epsii$ et si ce n'est pas un carr{\'e}, alors\/
$\thetabar(\Omega)$ engendre le groupe\/ $\zmodiixzmodii$.
\endth
\demonstration On a $\thetabar(\epsi)=(1,\machin)$.
D'autre part, le lemme~6 fournit un $\xi\in\ketoile$ tel que
$\xi^2-\epsi$ soit non~carr{\'e} ; on a alors $\xi^2\in\Omega$
(cf.~lemme~5) et $\thetabar(\xi^2)=(0,1)$.\cqfd

\th LEMME 12
\enonce
Supposons maintenant que\/ $\epsi\neq\epsii$.  Alors\/
$\thetabar(\Omega)$ engendre le groupe\/ $\zmodiixzmodii$.
\endth
\demonstration
Il y a quatre cas {\`a} consid{\'e}rer.

\puce $\epsi,\epsii$ sont tous deux des carr{\'e}s.  Alors le lemme~7,
appliqu{\'e} au couple $(\epsi,\epsii)$, fournit un $\xi_1\in\ketoile$
tel que $\alpha\xi_1^2\in\Omega$ et que
$\thetabar(\alpha\xi_1^2)=(1,0)$ ; le m{\^e}me lemme, appliqu{\'e} au
couple $(\epsii,\epsi)$, fournit un $\xi_2\in\ketoile$ tel que
$\alpha\xi_2^2\in\Omega$ et que $\thetabar(\alpha\xi_1^2)=(1,1)$

\puce $\epsi$ est un carr{\'e}, \ $\epsii$ ne l'est pas. On a
$\thetabar(0)=(1,\machin)$.  Si $\epsi-\epsii$ n'est pas un carr{\'e}
dans $\ketoile$, on a $\thetabar(\epsi)=(0,1)$. Sinon, le lemme~8
fournit un $\xi\in\ketoile$ tel que $\xi^2\in\Omega$ et que
$\thetabar(\xi^2)=(0,1)$.

\puce $\epsi$ n'est pas un carr{\'e}, \ $\epsii$ en est un. On a
$\thetabar(0)=(1,\machin)$. Si $\epsii-\epsi$ n'est pas un carr{\'e}
dans $\ketoile$, on a $\thetabar(\epsii)=(0,1)$. Sinon, le lemme~8
fournit un $\xi\in\ketoile$ tel que $\xi^2\in\Omega$ et que
$\thetabar(\xi^2)=(0,1)$.

\puce $\epsi,\epsii$ sont tous deux des non carr{\'e}s.  On a
$\thetabar(\epsi)=(1,\machin)$.  Si $-1$ est un carr{\'e} dans
$\ketoile$, on a $\thetabar(0)=(0,1)$.  Sinon, le lemme~9 fournit un
$\xi\in\ketoile$ tel que $\xi^2\in\Omega$ et que
$\thetabar(\xi^2)=(0,1)$.\cqfd
 
\goodbreak
\section{Le cas ramifi{\'e}}

Nous allons utiliser les r{\'e}sultats du \numero~7 avec
$\epsi=\eibar$, $\epsii=\eiibar$, $\eps_{1,2}=\eimeiibar$ dans la
d{\'e}finition~(11) de $\thetabar$ et $\Omega=\omega(M)$ (cf.\
\numero~6), ce qui est loisible gr{\^a}ce au lemme~4.

\goodbreak

\th LEMME 13 
\enonce
Les applications\/ $\theta$ et\/ $\thetabar\circ\omega$
co{\"\i}ncident sur $M\subset N$. 
\endth 
\demonstration C'est clair d'apr{\`e}s le lemme~3 et les
d{\'e}finitions.\cqfd 

R{\'e}sumons d'abord les r{\'e}sultats du calcul de
$\thetabar(\Omega)$ dans les trois cas~:

\rem Le cas i) $\ei/\eii\equiv1\,\modm$ et $\ei\in\normlk$\endrem  
Dans ce cas on a $\epsi=\epsii$ et cet {\'e}l{\'e}ment est un
carr{\'e} dans $\ketoile$.  Le lemme~10 affirme alors que
$\thetabar(\Omega)$ engendre le sous-groupe $\{0\}\times\zmodii$.

\rem Le cas ii) $\ei/\eii\equiv1\,\modm$ et $\ei\notin\normlk$\endrem Ici
$\epsi=\epsii$, mais ce n'est pas un~carr{\'e}. Le lemme~11 affirme
que $\thetabar(\Omega)$ contient $(0,1)$ et $(1,\machin)$.

\rem Le cas iii) $\ei/\eii\not\equiv1\,\modm$\endrem Ici $\epsi$ et
$\epsii$ sont distincts.  Il a {\'e}t{\'e} d{\'e}montr{\'e} au
lemme~12 que $\thetabar(\Omega)$ engendre le groupe $\zmodiixzmodii$.

Pour terminer la d{\'e}monstration de la prop.~2, il suffit de
remarquer qu'un {\'e}l{\'e}ment $x\in N$ qui n'est pas dans $M$ (donc
tel que $\xbar=\epsi$) a pour image $\theta(x)=(0,\machin)$ dans le
cas~{\it i\/}).\cqfd

\remarque\ Les courbes de genre~1 qui interviennent dans la
d{\'e}monstration du cas~{\it iii\/}) admettent un morphisme vers la
courbe ab{\'e}lienne~(5) qui s'{\'e}crit $x=\alpha\xi^2$,
$y=\alpha\xi\eta\zeta$ si $\epsi,\epsii$ sont tous deux des carr{\'e}s
et $x=\xi^2$, \ $y=\alpha\xi\eta\zeta$ sinon (avec les notations des
lemmes~7--9.)
 
\remarque\ C'est le groupe $\azeroo{X}$ qui est un invariant
birationnel de $X$ (\numero~3, th.~1) et non le plongement
$\azeroo{X}\fleche\zmodiixzmodii$.  En effet, soient
$\ei^\prime,\eii^\prime\in\Ketoile$ des {\'e}l{\'e}ments tels que
$\omega(-\eiprim)$ soit un carr{\'e} et $v(\eiprim)<v(\eiiprim)$.
Pour la surface de Ch{\^a}telet $Y$ (cf.~(2)) :
$$
y^2-dz^2=x\xprim(x-\eiprim\xprim)(x-\eiiprim\xprim)t^2,
$$
l'image de $\azeroo{Y}$ est {\'e}gale {\`a} $\zmodii\times\{0\}$ ; par
contre, en se ramenant au cas~{\it i\/}) de la prop.~2 par le
changement de variable $x_1=x-\eiprim$, on trouve
$\{0\}\times\zmodii$ comme image de $\azeroo{X}$, \ $X$ {\'e}tant la
surface de Ch{\^a}telet~(2) avec $\ei=-\eiprim$, \
$\eii=\eiiprim-\eiprim$.

\section{L'homomorphisme de restriction}

Dans ce \numero, nous d{\'e}montrons la prop.~3, par une division
assez naturelle en quatre cas suivant le type de ramification de $L$
et de $\Kprim$.  Comme dans les d{\'e}monstrations des propositions~1
et~2, on peut supposer que $X$ est une surface de Ch{\^a}telet,
donn{\'e}e par l'{\'e}quation~(2).

\goodbreak

Si $d$ est un carr{\'e} dans $\Kprim$, la surface
$\Xprim=X\times_K\Kprim$ est $\Kprim$-birationnelle {\`a} $\P_2$ et
$\azeroo{\Xprim}=0$ \citer\blochlivre( p.~7.1).  Comme le degr{\'e}
$[\Kprim:K]$ est alors pair, la prop.~3 dit que l'homomorphisme
$\azeroo{X}\fleche\azeroo{\Xprim}$ est nul, ce qui est trivialement
vrai.  On va supposer donc que $d$ n'est pas un carr{\'e} dans
$\Kprim$ non plus.

L'extension quadratique $\Lprim=\Kprim(\sqrt{d}\,)$ de $\Kprim$ est
ramifi{\'e}e si et seulement si $\vprim(d)$ est impair, o{\`u}
$\vprim$ est la valuation normalis{\'e}e de~$\Kprim$.  Il lui
correspond un homomorphisme surjectif
$\hprim:\Kprimetoile\fleche\zmodii$ de noyau $\normlprimkprim$.  Comme
nous allons le voir, tout revient~{\`a} la d{\'e}termination de la
restriction~$\hprim\,|_{\Ketoile}$.

Relativement {\`a} $\Xprim$ sur $\Kprim$, on d{\'e}signe par
$\Nprim\subset\Kprim$ l'analogue de l'ensemble $N$ ({\it
Remarque\/}~3) pour $X$ sur $K$ et par
$\thetaprim:\Xprim(\Kprim)\fleche\zmodiixzmodii$ l'analogue de
l'application~(9), ainsi que l'application qu'il induit sur
$\Nprim$.                 
                          
Abr{\'e}geons $\azeroo{X}$ en $A$ et $\azeroo{\Xprim}$ en $\Aprim$.
On dispose des applications $\theta:N\fleche A$ et
$\thetaprim:\Nprim\fleche\Aprim$ dont les images engendrent ({\it
Remarque\/}~3) ces groupes, et de l'inclusion $N\subset\Nprim$.  Il
s'agit de calculer l'image de $\thetaprim\,|_N$.

Il est clair qu'on peut supposer que l'extension $\Kprim\,|\,K$ est ou
bien non ramifi{\'e}e, ou bien totalement ramifi{\'e}e.

%\smallskipamount=2.0pt plus 1.0pt minus 1.0pt

\smallbreak
\puce{\it Supposons que les extensions\/ $L$ et\/ $\Kprim$ de\/~$K$ sont non
ramifi{\'e}es.}
\demonstration  Dans ce cas, le degr{\'e} $[\Kprim:K]$ ne peut pas {\^e}tre
pair puisque nous avons suppos{\'e} que $\Kprim$ est lin{\'e}airement
disjointe de $L$.  L'extension $\Lprim$ de $\Kprim$ est non
ramifi{\'e}e et comme $\vprim\,|_{\Ketoile}=v$, on a
$\thetaprim\,|_N=\theta$. Par suite, l'application
$A\fleche\Aprim$ est un isomorphisme.\cqfd

\remarque\ D'apr{\`e}s la prop.~1, le groupe $\azeroo{\Xprim}$
``~reste le m{\^e}me~'' que $\azeroo{X}$ lorsqu'on passe de $K$ {\`a}
une extension non ramifi{\'e}e $\Kprim$, qu'elle soit de degr{\'e}
pair ou impair.

\smallbreak
\puce{\it  Supposons que l'extension\/ $L$ est non ramifi{\'e}e et que\/
$\Kprim$ est ramifi{\'e}e.}\par
\demonstration  L'extension $\Lprim$ de $\Kprim$ est non ramifi{\'e}e. Comme
la restriction de $\thetaprim$ {\`a} $N$ est $\theta$ si $[\Kprim:K]$
est impair (puisqu'alors $v^\prime(x)$ a la m{\^e}me parit{\'e} que
$v(x)$, quel que soit $x\in\Ketoile$) et $0$ sinon (car $v^\prime(x)$
est pair quel que soit $x\in\Ketoile$), la prop.~3 est
d{\'e}montr{\'e}e dans ce cas.\cqfd

\remarque\  Noter que, d'apr{\`e}s la prop.~1, $\azeroo{X}$
``~reste le m{\^e}me~'' quand on passe de $K$ {\`a} $\Kprim$ sauf si
on est dans le cas~{\it iii\/}) (de la prop.~1) au d{\'e}part et si
$[\Kprim:K]$ est pair, auquel cas $\azeroo{X}$ est chang{\'e} de
$\zmodiixzmodii$ en $0$ si $v(\ei-\eii)=r$ (bonne r{\'e}duction
potentielle) et en $\zmodii$ si $v(\ei-\eii)>r$ (r{\'e}duction
multiplicative potentielle).

\goodbreak\smallbreak
\puce{\it  Supposons que l'extension\/ $L$ est ramifi{\'e}e et que\/
$\Kprim$ ne l'est pas.}\par
\demonstration 
L'extension quadratique $\Lprim$ de $\Kprim$ est ramifi{\'e}e. Soit
$\pi$ une uniformisante de $K$ adapt{\'e}e {\`a} l'extension $L$ ;
$\pi$ est alors une uniformisante de $\Kprim$ et elle est adapt{\'e}e
{\`a} $\Lprim$ (comme $\pi\in\normlk$, il appartient
aussi {\`a} $N_{L^\prime|\Kprim}(L^{\prime\times})$).  

On note $\omegaprim:\Kprimetoile\fleche\kprimetoile$ l'homomorphisme
attach{\'e} {\`a} l'uniformisante $\pi$ (cf.~\numero~5).  On associe
{\`a} $\epsi=\omegaprim(\ei)$, \ $\epsii=\omegaprim(\eii)$ et
$\eps_{1,2}=\omegaprim(\ei-\eii)$ une application
$\thetabarprim:\kprim\fleche\zmodiixzmodii$ de la m{\^e}me mani{\`e}re
que $\thetabar$ l'avait {\'e}t{\'e} (\numero~7) {\`a}
$\epsi=\omega(\ei)$, \ $\epsii=\omega(\eii)$ et
$\eps_{1,2}=\omega(\ei-\eii)$.  
\goodbreak

Notons $\Omega=\omega(M)$ et $\Omegaprim=\omegaprim(\Mprim)$, o{\`u}
$\Mprim\subset\Nprim$ contient pr{\'e}cis{\'e}ment $\ei$ et les
$x\in\Nprim$ tels que $\omegaprim(x)\neq\eibar$~(cf.~\numeros~6 et~8).
On a l'inclusion $\Omega\subset\Omegaprim$ et les diagrammes
$$\def\\{\mskip-2\thickmuskip} \def\'{\phantom{'}}
\diagram{
\Mprim &{\!\\\droite{\thetaprim}\\}&\Aprim &\phantom{\droite{}}
&\Mprim&{\!\\\droite{\omegaprim}\\}&\Omegaprim
&{\!\\\droite{\thetabarprim}\\}&\Aprim\cr
\dessus{\'}&&\dessus{\'}&&\dessus{\'}&&\dessus{\'}\cr
M\'&{\!\\\droite{\theta}\\}&A\'&&M\'&{\!\\\droite{\omega}\\}&\Omega\'\cr}
\abovedisplayskip=5.0pt plus 3.0pt minus 3.0pt
\belowdisplayskip=5.0pt plus 3.0pt minus 3.0pt
$$
sont commutatifs.  L'image de $M\fleche\Aprim$ se calcule donc via
$\thetabarprim\circ\omega$.  Mais la restriction de $\thetabarprim$
{\`a} $\Omega$ est {\'e}gale {\`a} 0 si $[\Kprim:K]$ est pair (car
$\ketoile\subset k^{\prime\times2}$) et {\`a} $\thetabar$ sinon (car
$k^{\prime\times2}\cap\ketoile=k^{\times2}$), d'o{\`u} la prop.~3 dans
ce cas.\cqfd

\remarque\  Noter que, d'apr{\`e}s la prop.~2, $\azeroo{X}$
``~reste le m{\^e}me~'' quand on passe de $K$ {\`a} $\Kprim$ sauf si on est
dans le cas~{\it ii\/}) (de la prop.~2) au d{\'e}part et si
$[\Kprim:K]$ est pair, auquel cas il est chang{\'e} de
$\zmodiixzmodii$ en $\zmodii$ (potentiellement r{\'e}duction
multiplicative d{\'e}ploy{\'e}e).

\smallbreak
\puce{\it Supposons que les extensions\/ $L$ et\/ $\Kprim$ sont toutes deux
ramifi{\'e}es.}\par
\demonstration 
Consid{\'e}rons d'abord le cas o{\`u} $n=[\Kprim:K]$ est impair.
L'extension $\Lprim$ de $\Kprim$ est alors ramifi{\'e}e.  La
restriction de $\hprim$ {\`a} $\Ketoile$ est $h$~: il suffit de faire
voir que $\hprim\,|_{\Ketoile}$ n'est pas nulle.  Prenons une
uniformisante $\varpi$ de $K$ qui n'est pas adapt{\'e}e {\`a} $L$,
c'est-{\`a}-dire telle que $\varpi\notin\normlk$.  Alors $\varpi\notin
N_{L^\prime|\Kprim}(L^{\prime\times})$, car sinon $\varpi^n=
N_{\Kprim|K}(\varpi)$ appartiendrait {\`a}
$N_{L^\prime|K}(L^{\prime\times})$ et {\it a fortiori\/} {\`a}
$\normlk$, mais $h(\varpi^n)=nh(\varpi)=1$.  Comme $\hprim(\varpi)=1$,
la restriction de $\hprim$ {\`a} $\Ketoile$ n'est pas nulle et l'on a
bien $\hprim\,|_{\Ketoile}=h$.  Par suite, l'application
$A\fleche\Aprim$ est surjective.

Si $n=[\Kprim:K]$ est pair, l'extension $\Lprim$ de $\Kprim$ est {\it
non ramifi{\'e}e}.  Comme $\vprim(x)$ est pair quel que soit
$x\in\Ketoile$, on a $\thetaprim(M)=0$ et l'application
$A\fleche\Aprim$ est nulle.\cqfd

\remarque\ Si l'extension $\Kprim\,|\,K$ est de degr{\'e} pair, le
groupe $\azeroo{X}$ ``~reste le m{\^e}me~'' quand on passe de $K$
{\`a} $\Kprim$ si on est dans le cas~{\it i\/}) (de la prop.~2) au
d{\'e}part.  Il est chang{\'e} de $\zmodiixzmodii$ en $0$ dans le
cas~{\it iii\/}) (bonne r{\'e}duction potentielle) et
en $\zmodii$ dans le cas~{\it ii\/}) (potentiellement
r{\'e}duction multiplicative).

\unvbox\bibbox
\bye